\documentclass{amsart}
\usepackage{graphicx}
\usepackage{amssymb}

\newtheorem{proposition}{Proposition}[section]
\newtheorem{theorem}{Theorem}[section]
\newtheorem{lemma}[theorem]{Lemma}
\newtheorem{cor}[theorem]{Corollary}

\theoremstyle{definition}

\theoremstyle{remark}

\numberwithin{equation}{section}

\newcommand{\A }{\mathcal{A}}
\newcommand{\Aql }{\overrightarrow{\A }}

\newcommand{\B }{\mathcal{B}}

\newcommand{\cA }{\mathcal{A}}

\newcommand{\cO}{\mathcal{O}}

\newcommand{\comp }{\mathbb{C}}
\newcommand{\C }{\mathbb{C}}
\newcommand{\dif }{\mathrm{d}}

\newcommand{\F}{\mathcal{F}}

\newcommand{\fsl}{\mathfrak{sl}}
\newcommand{\Gam }{\varGamma }
\newcommand{\gamgr }{\Gamma^1_q(\mathrm{Gr}(r,N))}

\newcommand{\kopr }{\varDelta }
\newcommand{\kow }{\varDelta }
\newcommand{\K }{K}
\newcommand{\Kkow }{\varDelta _\K }
\newcommand{\lact }{\triangleright}
\newcommand{\lid }{\mathcal{L}}
\newcommand{\Lin }{\mathrm{Lin}}

\newcommand{\N}{\mathbb{N}}

\newcommand{\ot }{\otimes }

\newcommand{\oqg}{{\mathcal O}_q(G)}

\newcommand{\oqgr }[1][r,N]{\mathcal{O}_q(\mathrm{Gr}(#1))}

\newcommand{\Oqmat }{\mathcal{O}_q(\mathrm{Mat}(r,N{-}r))}

\newcommand{\pair }[2]{\langle #1,#2\rangle }

\newcommand{\qd }{\hat{q}}
\newcommand{\qlm }[1]{\overrightarrow{#1}}

\newcommand{\R }{\mathbb{R}}
\newcommand{\real }{\mathbb{R}}

\newcommand{\SLN }{\mathcal{O}_q(\mathrm{SL}(N))}

\newcommand{\SUN }{\mathcal{O}_q(\mathrm{SU}(N))}

\newcommand{\U }{U}
\newcommand{\ug }{U_q(\mathfrak{g})}

\newcommand{\UslN }{U_q(\mathfrak{sl}_N)}
\newcommand{\Ubar }{\overline{\U }}
\newcommand{\vep }{\varepsilon }

\newcommand{\Z }{\mathbb{Z} }

\begin{document}

\title[Differential Calculus on Quantum Grassmannians II]
 {Differential Calculus on Quantum Complex Grassmann Manifolds II:
   Classification}

\author{Istv\'an Heckenberger}
\address{Mathematisches Institut, Universit\"at Leipzig, Augustusplatz 10,
         04109 Leipzig, Germany}

\email{heckenbe@mathematik.uni-leipzig.de} 

\author{Stefan Kolb}
\email{kolb@itp.uni-leipzig.de}

\thanks{S.~K.  was supported by the Deutsche Forschungsgemeinschaft
  within the scope of the postgraduate scholarship programme
  ``Graduiertenkolleg Quantenfeldtheorie'' at the
  University of Leipzig}



\begin{abstract}
For differential calculi over certain right coideal subalgebras of
quantum groups the notion of quantum tangent space is introduced.
In generalization of a result by Woronowicz a one to one correspondence
between quantum tangent spaces
and covariant first order differential calculi is established. This result is
used to classify differential calculi over quantum Grassmann manifolds $\oqgr$.
It turns out that up to special cases in low dimensions there exists
exactly one such calculus of classical dimension $2r(N{-}r)$.
\end{abstract}

\maketitle

In the framework of quantum groups and quantum spaces there appear many
examples of $q$-deformed coordinate algebras which allow well behaved
covariant first order differential calculi (FODC) in the sense of Woronowicz
\cite{a-Woro2},\cite{b-CP94},\cite{b-KS}. On the other hand there exists
no general construction of a deformation of classical K\"ahler differentials
in this setting. Therefore the task of classification of all covariant FODC
over quantum spaces naturally arises and has been settled for many examples
\cite{a-Po92},\cite{a-ApelSchm94},\cite{a-SchSch1},\cite{a-SchSch2},
\cite{a-BS98},\cite{a-HeckSchm2},\cite{a-Maj98},\cite{a-Welk98},
\cite{a-herm01}.

There are several techniques to perform classification. In
\cite{a-PuWo89} W.~Pusz and S.L.~Woronowicz consider calculi over
quantum vector spaces generated as left modules by the differentials of the
generators. A general ansatz is made and coefficients are determined by
covariance. This approach has been applied to other quantum spaces
\cite{a-Po92}, \cite{a-ApelSchm94}, \cite{a-Welk98},
yet the calculational effort of this method becomes very large for more
involved examples.

For any Hopf algebra $\cA$ there exists a one to
one correspondence between differential calculi and certain right (or left)
ideals of $\cA^+=\ker(\vep)\subset\cA$ \cite{a-Woro2}.
Right ideals have been used in \cite{a-SchSch1}, \cite{a-SchSch2},
\cite{a-BS98} to classify bicovariant differential calculi over quantum groups.
In \cite{a-herm01} U.~Hermisson has generalized this
method to certain right (or left) coideal subalgebras $\B\subset\cA$ and
classified all 2-dimensional covariant FODC over Podle\'s' quantum sphere.

A reformulation of the method of right ideals involves the notion of
quantum tangent space. In the case of a Hopf algebra $\cA$ the quantum
tangent space $T_\Gamma$ of a covariant FODC $\Gamma$ is a subspace of the
dual Hopf algebra $\cA^\circ$ which determines the calculus $\Gamma$
uniquely. Quantum tangent spaces have been used for
classification in \cite{a-HeckSchm2}, \cite{a-Maj98}.
The advantage of the quantum tangent space approach in the case of coordinate
algebras of quantum groups $\cA=\oqg$ stems from the fact that $\cA^\circ$
possesses a well understood subalgebra $\Ubar$ which is closely related to
the $q$-deformed universal enveloping algebra $U=\ug$. The main strategy is
to reduce the classification problem of covariant FODC to a classification
problem of quantum tangent spaces in $\Ubar$.

In this paper the quantum tangent space method is generalized to quantum
spaces. More precisely we consider right coideal subalgebras
$\B$ of a quantum group $\cA=\oqg$ obtained as right $K$ invariants for
certain left coideal subalgebras $K$ of $U$. In this situation a quantum
tangent space is a subspace of the dual coalgebra $\B^\circ$.
Now the strategy for classification is the following.
One has to find a suitable sufficiently small subcoalgebra
$\Ubar\subset\B^\circ$. Using representation theory of $\B$ one has to
show that any quantum tangent space lies in $\Ubar$.
Finally quantum tangent spaces in $\Ubar$ have to be classified.
The crucial point of this strategy is the right choice of $\Ubar$.

Here this program is performed for quantum Grassmann
manifolds $\oqgr$ \cite{a-NDS97}, \cite{a-DS99}, \cite{a-Kolb} in the so
called quantum subgroup case. All covariant FODC over $\oqgr$ of dimension
up to $2r(N{-}r)$ are classified.
It turns out that up to special cases in low dimensions there exists
exactly one covariant FODC which has the same dimension as its classical
counterpart. This calculus has been constructed and investigated
in \cite{a-Kolb}.

The ordering of this paper is as follows. In Section \ref{tangent} the 
notion of quantum tangent space for a certain class of quantum spaces is
introduced. A one to one correspondence between covariant FODC and
quantum tangent spaces is established. This result strongly relies on the
identification of covariant FODC with certain left ideals of $\B^+$ given
in \cite{a-herm01}.

Quantum Grassmann manifolds $\oqgr$ are recalled in Section \ref{q-grass}.

In Section \ref{uk+u} the structure of $\Ubar=U/K^+U$ in the case
$\B=\oqgr$ is investigated. It is shown that there exists a nondegenerate
pairing between $\B/(\B^+)^{k+1}$ and the set $\Ubar_k\subset\Ubar$ of elements
of degree $k$ with respect to the coradical filtration of $\Ubar$.

Section \ref{gradrep} is devoted to the representation theory of $\oqgr$.
Technical lemmata are obtained by explicit calculations using generators and
relations of $\oqgr$. For this reason some of the relations are collected in
Appendix \ref{relations}.

Combining these results it is shown in Section \ref{classification} that any
element of a finite dimensional quantum tangent space vanishes on $(\oqgr^+)^k$
for some $k$ and therefore belongs to $\Ubar$. Finally all quantum tangent
spaces in $\Ubar$ of dimension up to $2r(N{-}r)$ are determined

If not stated otherwise all notations and conventions coincide with those
introduced in \cite{b-KS}.

\section{Quantum tangent space}\label{tangent}

Let $\U $ denote a Hopf algebra with bijective antipode and $\K \subset \U $
a left coideal subalgebra, i.e.
$\Kkow :\K \rightarrow \U \ot \K $. Consider a
tensor category $\mathcal{C}$ of finite dimensional left $\U $-modules. Let
$\A :=\U ^0_{\mathcal{C}}$ denote the dual Hopf algebra generated by the
matrix coefficients of all $\U $-modules in $\mathcal{C}$.
Assume that $\A$ separates the elements of $\U$.
Define a right coideal subalgebra $\B \subset \A $ by
\begin{align}\label{Bdef}
  \B :=\{b\in \A \,|\,\pair{k}{b_{(1)}}b_{(2)}=0 \quad
  \text{for all $k\in \K ^+$}\},
\end{align}
where $K^+=\{k-\vep(k)\,|\,k\in K\}$.
Assume $\K $ to be $\mathcal{C}$-semisimple, i.e. the restriction of any
$\U$-module in $\mathcal{C}$ to the subalgebra $K\subset \U$ is isomorphic
to the direct sum of irreducible $K$-modules. By \cite{a-MullSch}
Theorem 2.2 (2) this implies that $\A $ is a faithfully flat $\B $ module.
Let $\Aql $ denote the left $\A $-module coalgebra $\Aql :=\A /\A \B ^+$,
where $\B^+=\{b\in \B\,|\,\vep(b)=0\}$.
Further for any right $\B $-module $M$ set $\qlm{M}:=M/M\B ^+$, and for any
right $\A $-comodule $\kow:N\rightarrow N\otimes \A $ define
$\qlm{\kow}:N\rightarrow N\otimes \Aql$ to be the compositon of $\kow$ with
the canonical projection.

\begin{lemma}\label{l-nondeg}
The pairing $\K \times \A /\A \B ^+ \rightarrow \comp $ is nondegenerate.
\end{lemma}

\begin{proof}
It is shown in the proof of \cite{a-MullSch} Theorem 2.2 (1), (2) that $\A /\A \B ^+$
is equal to the image of $\A $ under the restriction map
$\U ^0\rightarrow \K ^0$.
\end{proof} 

Let $\lid \subset \B ^+$ be a subspace of finite codimension. Then the
orthogonal complement $\lid $ in the dual vector space $\B'$ of $\B$ is
defined by
\begin{align}
  T_\lid :=\{f\in \B '\,|\,f(x)=0\quad \text{for all $x\in \lid $}\}.
\end{align}
Obviously, $\vep \in T_\lid $. If $\lid $ is the left ideal corresponding
to a unique right-covariant first order differential calculus $\Gam $
over $\B $, see \cite{a-herm01}, then $T_\lid ^+=\{t\in T_\lid\,|\,t(1)=0\}$ is
called the \textit{quantum tangent space} of $\Gam $.

On the other hand, for any finite dimensional subspace $T\subset \B '$,
$\vep \in T$, define
\begin{align}
\lid _T:=\{b\in\B \,|\,f(b)=0\quad \text{for all $f\in T$}\}.
\end{align}
Then $\lid _T\subset \B ^+$.

\begin{proposition} Let $\lid $ and $T_\lid $ be given as above, then
\begin{enumerate}
\item $(\B \lid \subset \lid \mbox{ and } T_\lid \subset \B ^\circ )
       \Rightarrow \kopr T_\lid \subset \B ^\circ \ot T_\lid $,
\item $\qlm{\kow } \lid \subset \lid \ot \Aql \Rightarrow
      T_\lid K\subset T_\lid ,$
\item $(\B \lid \subset \lid \mbox{ and }\qlm{\kopr } \lid \subset
       \lid \ot \Aql ) \Rightarrow T_\lid \subset \B ^\circ ,$
\item $\dim _\comp T_\lid =\dim _\comp \B ^+/\lid +1$,
\item $\lid _{T_\lid }=\lid $.
\end{enumerate}
\end{proposition}

\begin{proof}
1. Let $f\in T_\lid \subset \B ^\circ $ and consider
$\kopr f=\sum _if_i\ot g_i$, where $f_i$ are linearly independent in
$\B ^\circ $. By assumption $0=f(bx)=\sum _if_i(b)g_i(x)$
for all $b\in \B ,x\in \lid $.
Moreover, there exist $b_j\in B$ such that $f_i(b_j)=\delta _{ij}$ for all
$i,j$. Thus $g_i(x)=0$ for all $x\in \lid $.

2. Assume $k\in \K ,f\in T_\lid $ and $x\in \lid $. Since the pairing of $\K $
and $\Aql $ is well-defined, we get
$(fk)(x)=f(x_{(0)})k(x_{(1)})=0$ as $f(\lid )=0$.

4. and 5. This follows from the fact that the codimension of
$\lid $ in $\B $ is finite.

3. By assumption and Theorem 2 of \cite{a-herm01}, and exchanging left and
right, $\lid $ corresponds to a right-covariant
FODC $\Gam $ over $\B $. Choose a basis $\dif x_1,\dots, \dif x_k$ of
$\qlm{\Gam}=\Gam /\Gam \B ^+$ and define functionals
$\chi ^i,i=1,\ldots ,k$, in $\B '$ by
\begin{align}
\overline{\dif b}=\dif x_i\chi ^i(b).
\end{align}
Here $\overline{\dif b}$, $b\in \B $, denotes the element of
$\qlm{\Gam}$ represented by $\dif b$.
Since $\Gam =\Lin \{\dif a\,b\,|\,a,b\in \B \}$, we get $\qlm{\Gam}
=\Lin \{\overline{\dif b}\,|\,b\in \B \}$. Moreover, $\chi ^i(1)=0$
for all $i$. Therefore the functionals
$\{\chi ^1,\ldots ,\chi ^k,\vep \}$
are linearly independent.

Let $y_i\in \B $ such that
$\overline{\dif y_i}=\dif x_i$. Then the left $\B $-module structure of
$\Gam $ induces a finite dimensional representation $\rho $ of $\B $ by
\begin{align}
\overline{b\dif y_i}=\dif x_j\rho^j_i(b).
\end{align}
By the Leibniz rule
\begin{align}
\chi ^i(ab)=\chi ^i(a)\vep (b)+\rho ^i_j(a)\chi ^j(b)
\end{align}
which implies $\chi ^i\in \B ^\circ $.
On the other hand in terms of $\Gam $
\begin{align}
\lid =\{a_i^+\vep (b_i)\,|\,\dif a_i\,b_i=0\}.
\end{align}
Thus $\chi ^i(x)=0$ for all $x\in \lid $, i.\,e.\
$\chi ^i\in T_\lid $.
Now, $\dim _\comp T_\lid =\dim _\comp \B ^+/\lid +1
=\dim _\comp \Gam /\Gam \B ^+ +1$ by 4.
Therefore $\{\chi ^1,\ldots ,\chi ^k,\vep \}$ is a sufficiently large
set to span all of $T_\lid $.
\end{proof}  

\begin{proposition} Let $T$ and $\lid _T$ be given as above. Then 
\begin{enumerate}
\item $(T\subset \B ^\circ \mbox{ and }\kopr T\subset \B ^\circ \ot T)
       \Rightarrow \B \lid _T\subset \lid _T$,
\item $TK\subset T\Rightarrow \qlm{\kopr } \lid _T\subset
       \lid _T\ot \Aql $,
\item $\dim _\comp T=\dim _\comp \B ^+/\lid _T+1$,
\item $T_{\lid _T}=T$.
\end{enumerate}
\end{proposition}  

\begin{proof}
1. Let $x\in \lid_T , b\in \B $ and $f\in T$. Then
$f(bx)=f_{(1)}(b)f_{(2)}(x)=0$ by assumption.

2. Since $\B $ is a right $\A $ comodule, it is also a right $\Aql $
comodule with coaction $\qlm{\kopr }$. Let $x\in \lid_T $ and
$\qlm{\kopr }x=\sum _ix_i\ot y_i\in \B \ot \Aql $.
Assume $y_i$ to be linearly independent in $\Aql $. Then by
Lemma \ref{l-nondeg} there exist $k_j\in \K $ such that
$k_j(y_i)=\delta _{ij}$.
Now for any $t\in T$ one has $t(x_j)=(tk_j)(x)=0$ by assumption.

3. The subspace of $\B $ where all elements of $T$ vanish has
codimension $\dim _\comp T$.

4. This claim follows for instance from $T\subset T_{\lid _T}$ and
$\dim _\comp T=\dim _\comp T_{\lid _T}$.
\end{proof}

\begin{cor}\label{corresp}
Under the above assumptions there is a one to one correspondence between
$n$-dimensional covariant FODC over $\B $ and $(n+1)$-dimensional
subsets $T^\vep\subset \B ^\circ $ such that
\begin{align}
\vep \in T^\vep,\quad \Delta T^\vep \subset \B ^\circ \ot T^\vep,\quad T^\vep
 K\subset T^\vep.
\end{align}
\end{cor}

\section{Quantum Grassmann manifolds}\label{q-grass}
In the remaining sections of this paper all considerations will be restricted
to the following example of quantum complex Grassmann manifolds.
Let $q\in \real $ be transcendental. This assumption is needed only in the
proofs of Lemma \ref{EFprim} and of the main Theorem \ref{mainTh} where
duality between $\Ubar_+$ and $\Oqmat$ is used. All other argumentations
only use $q\in\R\setminus\{-1,0,1\}$.
Consider the $q$-deformed universal enveloping algebra $\U:=\UslN $.
The subalgebra $\K $ generated by
\begin{align}
\{E_i,F_i,K_j\,|\,i\neq r,j=1,\dots,N-1\}\label{Kdef}
\end{align}
is a coalgebra and hence a left coideal. The category $\mathcal{C}$ of type
one representations of $\UslN $ is a tensor category. The matrix coefficients
of $\mathcal{C}$ generate $\SLN $ and $\K $ is $\mathcal{C}$-semisimple.
Moreover, since $q$ is not a root of unity, the pairing between
$\SLN$ and $\UslN$ is nondegenerate.
The subalgebra $\B $ defined by (\ref{Bdef}) is just the $q$-deformed
coordinate algebra $\oqgr $ of the Grassmann manifold of $r$-dimensional
subspaces in $\comp ^N$ \cite{a-NDS97}, \cite{a-DS99} in the quantum subgroup
case, i.e.~the case when $K$ is a Hopf subalgebra of $\U$.

It has been shown in \cite{a-Kolb} that $\oqgr$ can be written in terms of
generators $z_{ij}$, $i,j=1,\dots, N$ which satisfy the relations
\begin{align}
q^{2b-2j}z_{ij}z_{ka}\hat{R}^{kb}_{jl}\hat{R}^{il}_{cd}&=q^{2l-2i}z_{ci}z_{jk}
     \hat{R}^{jl}_{id}\hat{R}^{ab}_{kl} \label{refl}\\
\sum_{i=1}^N z_{ii}&=\frac{1-q^{-2s}}{q-q^{-1}}\label{trace}\\
       q^{2N+1}\sum_{n=1}^Nq^{-2n}z_{in}z_{nk}&=z_{ik},\label{proj}
\end{align}
where $s=N{-}r$.
The explicit form of some of the relations (\ref{refl}) -- (\ref{proj}) is
given in Appendix \ref{relations} and will be used in Section \ref{gradrep}.

The right coaction of $\SLN$ on the generators $z_{ij}$ is given by
\begin{align}
  \kow z_{ij}=z_{kl}\ot u^k_iS(u^j_l)\label{zkow}
\end{align}
where $u^i_j$, $i,j=1,\dots,N$ denote the matrix elements of the
fundamental corepresentation
of $\SLN$. One obtains an induced left action of $\UslN$ which in the
convetions of \cite{b-KS}, Sect.~8.4.1 reads as
\begin{align}
E_k\lact z_{ij}&=\delta_{ik}q^{\delta_{kj}-\delta_{k,j-1}}z_{i+1,j}-
                 \delta_{j,k+1}qz_{i,j-1}\label{Eact}\\
F_k\lact z_{ij}&=\delta_{i,k+1}z_{i-1,j}-
                 \delta_{j,k}q^{\delta_{ik}-\delta_{i,k+1}-1} z_{i,j+1}\\
K_k\lact z_{ij}&=q^{\delta_{jk}-\delta_{j,k+1}-\delta_{ik}+\delta_{i,k+1}}
                 z_{ij}.
\end{align}
The embedding $i:\oqgr\rightarrow \SLN$ is given by
\begin{align*}
  i(z_{ij})=\sum_{i=r+1}^Nq^{-2N-1+2k}u^k_i S(u^j_k),
\end{align*}  
thus in particular
\begin{align}\label{zeps}
  \vep(z_{ij})=\begin{cases}
                 q^{-2N-1+2i} & \text{if $i=j>r$,}\\
                 0 & \text{else.}
               \end{cases}  
\end{align}

\section{The coalgebra $U/K ^+U$  }\label{uk+u}

Recall from \cite{a-DCoKa90} that there exists a filtration of $U$ such that
the associated graded algebra is $q$-commutative, i.e. given by generators
$E_\alpha,F_\beta,K^\pm_j$ and relations $t_1t_2=q^n t_2t_1$ for all
$t_1,t_2\in\{E_\alpha,F_\beta,K^\pm_j\}$ and some $n=n(t_1,t_2)\in\Z$ .
Here $\alpha, \beta\in \Phi^+$ denote the positive roots of $\fsl_N$.
Fix a reduced decomposition of the longest element $w$ of the Weyl group.

\begin{lemma}\label{baslem}
 Let $g_i,i=1,2,\dots$ denote the generators $E_\alpha,F_\beta,K_j$ of $U$
 with respect to $w$ in an arbitrary order. Then the elements
 \begin{equation*}
   \prod_i g_i^{n_i}
 \end{equation*}  
 $n_i\in \N_0$ if $g_i=E_\alpha,F_\alpha$ and $n_i\in\Z$ if $g_i=K_j$ form
 a vector space basis of $U$.
\end{lemma}
\begin{proof}
  This follows from the $q$-commutativity \cite{a-DCoKa90}, Prop 1.7 d.
\end{proof}  
Write the set $\Phi^+$ of positive roots as a disjoint union
\begin{equation*}
  \Phi^+=\Phi^+_r\cup \Phi^+_{-r}
\end{equation*}
where $\Phi^+_r=\{\alpha_{ij}=\sum_{k=i}^j\alpha_k|i\le r\le j\}$ and
$\alpha_k$ denote the simple roots.
\begin{proposition}\label{basprop}
  Let $\beta_i$ and $\beta_i'$ denote the elements of $\Phi^+_r$ in
  arbitrary fixed orders. Then the elements
  $$\prod_i (F_{\beta_i})^{n_i}\prod_j (E_{\beta_j'})^{m_j}$$
  $n_i,m_j\in\N_0$ form a vector space basis of $U/K^+U$.
\end{proposition}
\begin{proof}
  Let $\gamma_i$ and $\gamma_i'$ denote the elements of $\Phi^+_{-r}$ in
  arbitrary fixed orders. By Lemma \ref{baslem} the elements
  \begin{equation}
    K_1^{i_1}\dots K_{N-1}^{i_{N-1}}\prod_i (F_{\gamma_i})^{r_i}
    \prod_j (E_{\gamma_j'})^{s_j}\prod_k (F_{\beta_k})^{n_k}
    \prod_l (E_{\beta_l'})^{m_l}\label{basis}
  \end{equation}
  $i_k\in\Z,r_i,s_j,n_k,m_l\in\N_0$ form a vector space basis of $U$.
  Thus it suffices to show that the elements
  \begin{eqnarray}
     (K_1^{i_1}\dots K_{N-1}^{i_{N-1}}-1)\prod_k (F_{\beta_k})^{n_k}
    \prod_l (E_{\beta_l'})^{m_l}\label{bas1}\\
     K_1^{i_1}\dots K_{N-1}^{i_{N-1}}\prod_i (F_{\gamma_i})^{r_i}
    \prod_j (E_{\gamma_j'})^{s_j}\prod_k (F_{\beta_k})^{n_k}
    \prod_l (E_{\beta_l'})^{m_l}\label{bas2}
  \end{eqnarray}
   $i_k\in\Z,r_i,s_j,n_k,m_l\in\N_0$, $\sum r_i+\sum s_j\ge 1$, form a
   vector space basis of $K^+U$.
   Indeed the expressions (\ref{bas1}) and (\ref{bas2}) form a set of
   linearly independent elements of $K^+U$. Any element of $K^+U$ can be
   written as a sum of expressions of the form
    \begin{eqnarray}
    G K_1^{i_1}\dots K_{N-1}^{i_{N-1}}\prod_i (F_{\gamma_i})^{r_i}
    \prod_j (E_{\gamma_j'})^{s_j}\prod_k (F_{\beta_k})^{n_k}
    \prod_l (E_{\beta_l'})^{m_l}\label{Gexpression}
  \end{eqnarray}
  where $G\in\{K_i-1,F_j,E_j|j\neq r\}$ and
  $i_k\in\Z,r_i,s_j,n_k,m_l\in\N_0$.
  If $G=E_j$ then $G$ $q$-commutes with $K_1^{i_1}\dots K_{N-1}^{i_{N-1}}$.
  If $\prod_i(F_{\gamma_i})^{r_i}\neq F_j$ then by
  reordering $E_j\prod_i(F_{\gamma_i})^{r_i}$ according to the above basis
  (\ref{basis}) one obtains monomials of the form
  $$ K_1^{i_1'}\dots K_{N-1}^{i_{N-1}'}\prod_i(F_{\gamma_i})^{r_i'}E_j^\delta$$
  where $\delta=1$ or $\sum r_i'\neq 0$. As the elements $E_{\gamma_j'}$
  generate a subalgebra with basis  $\prod_j (E_{\gamma_j'})^{s_j}$ and
  $E_j$ is an element of this subalgebra the expression (\ref{Gexpression})
  for $G=E_j$ can indeed be written as a linear combination of elements
  of the form (\ref{bas2}). If on the other hand
  $\prod_i(F_{\gamma_i})^{r_i}= F_j$ then the relation
  $$E_j F_j-F_j E_j=\frac{K_j-K_j^{-1}}{q-q^{-1}}$$
  implies the claim. The cases $G=K_i-1$ and $G=F_j$ are dealt with in a
  similar way.
\end{proof}

Let $\Ubar $ denote $\U /\K ^+\U $.
By Corollary 5.3.5 in \cite{b-Montg93} $\Ubar $ is pointed.
Recall that the coradical $U_0$ of $\UslN$ is the subalgebra generated
by the elements $K_i,i=1,\dots,N-1$, \cite{b-Montg93}, Lemma 5.5.5.
By \cite{b-Sweedler}, p.~182, Ex.~4 the coradical of $\Ubar$ is contained in
$\pi(U_0)$ where $\pi:\UslN \rightarrow \Ubar$ denotes the canonical
projection. Thus $\Ubar$ is connected, i. e. the coradical of $\Ubar$ is
equal to $\C \cdot 1$. 

For any coalgebra $C$ let  $P(C )=\{x\in C \,|\,\kopr x=1\ot x+x\ot 1\}$
denote the vector space of primitive elements of $C$.

\begin{lemma}\label{invK}
$P(\Ubar )K\subset P(\Ubar )$.
\end{lemma}

\begin{proof}
Since $K$ is a coalgebra, for $p\in P(\Ubar )$, $k\in K$ we get
\begin{align*}
\kopr (pk)&=pk_{(1)}\ot k_{(2)}+k_{(1)}\ot pk_{(2)}\\
&=pk_{(1)}\ot \vep (k_{(2)})1+\vep (k_{(1)})1\ot pk_{(2)}=pk\ot 1+1\ot pk.
\end{align*}
\end{proof}

Set $\Ubar _-=\Lin \{\prod _iF_{\beta _i}^{m_i}\}$,
$\Ubar _+=\Lin \{\prod _iE_{\beta '_i}^{n_i}\}$, where the products are
taken over all $i$ such that $\beta_i,\beta_i'\in\Phi_r^+$.
As  $\Ubar _+$ (resp. $\Ubar _-$) is the image of the Hopf subalgebra
$U_q(\mathfrak{b}^+)$ (resp. $U_q(\mathfrak{b}^-)$) under the canonical
projection $\pi$, the subspace $\Ubar _+\subset \Ubar$
(resp. $\Ubar _-\subset \Ubar$) is a subcoalgebra.

In what follows several $U$-module coalgebra filtrations of $\Ubar$ and
$U_q(\mathfrak{b}_\pm)$-module coalgebra filtrations of $\Ubar_{\pm}$ will
prove quite useful. Let $\mathcal{F}_1$ denote the filtration of $\Ubar$
defined by
\begin{align*}
  \deg_1\left( \prod_i (F_{\beta_i})^{n_i}\prod_j (E_{\beta_j'})^{m_j}\right)
  =\sum_in_i+ \sum_j m_j.
\end{align*}
The corresponding filtration of $U$ is defined by
\begin{align*}
  \deg _1(E_i)=\delta_{i,r}=\deg _1(F_i),\quad \deg_1(K_i)=0.
\end{align*}
This induces a filtration on $\Ubar_\pm$ which will also be denoted by
$\mathcal{F}_1$. 
Further let $\mathcal{F}_2$ denote the filtration of $\Ubar$ defined by
\begin{align*}
  \deg_2(E_i)=\deg_2(F_i)=1.
\end{align*}
The corresponding filtration of $U$ is defined by
\begin{align*}
  \deg _2(E_i)=1=\deg _2(F_i),\quad \deg_2(K_i)=0.
\end{align*}
Finally there exists a $\Z^{N-1}$ grading of $\Ubar$ and $\Ubar_\pm$ induced
by the standard $\Z^{N-1}$ grading of $\U$.

\begin{lemma}\label{EFprim}
  The primitive elements in $\Ubar_+$ and $\Ubar_-$ are given by
\begin{align*}
  P(\Ubar_+ )=\Lin \{E_{\beta '_i} \, | \, \beta_i'\in\Phi_r^+\},\quad
  P(\Ubar_- )=\Lin \{F_{\beta _i} \, | \, \beta_i\in\Phi_r^+\}.
\end{align*}  
\end{lemma}

\begin{proof} The $\U$-module coalgebra $\Ubar_{\pm}$ is graded with respect
  to $\mathcal{F}_1$. Let$(\Ubar_{\pm})_k$ denote the elements of degree $k$.
  As $(\Ubar_{\pm})_0=\C\cdot 1$ all elements of $(\Ubar_{\pm})_1$ have to
  be primitive.
  It remains to show that there are no primitive elements in the graded
  components $(\Ubar_{\pm})_k$ for $k>1$. The proof is carried out for the
  case $(\Ubar_{+})_k$.

  Let $V(0)_+$ denote the right highest weight $\U$-module given by one
  generator $v$ and relations
  \begin{align*}
    v F_i=0,\quad v K_i=v,\quad v E_j=0 \quad \text{for all $i,j=1,\dots, N-1,$
    $j\neq r.$ }
  \end{align*}  
  The module $V(0)_+$ can be endowed with a right $\U$-module coalgebra
  structure by (cf. \cite{a-SinVa}, \cite{a-SinShklyVa})
  \begin{align*}
    \kow v=v\otimes v.
  \end{align*}  
  Note that the coalgebras $\Ubar_+$ and $V(0)_+$ are both isomorphic to
  \begin{align*}
    U_q(\mathfrak{b}_+)\big/\langle(K_i-1) U_q(\mathfrak{b}_+),
                   E_j U_q(\mathfrak{b}_+)\,|\,j\neq r\rangle.
  \end{align*}  
It has been shown in \cite{a-SinShklyVa} that the $\U$-module coalgebra
$V(0)_+$ is
the graded dual of the $U$-module algebra of $q$-deformed functions on
$r\times (N{-}r)$ matrices $\Oqmat$.
Recall that the homogeneous component $(V(0)_+)_k$ of degree $k$ elements
is spanned by all monomials $\prod_i (E_{\beta_i'})^{n_i}$ with
$\sum_i n_i=k$, while the degree of monomials in $\Oqmat$ is given by the
number of factors.
Thus Lemma \ref{invK} implies that  $P(\Ubar)$ is a graded vector space.
Assume $x\in P(\Ubar_+)$ to be a homogeneous primitive element with
$\deg_1(x)>1$. As any monomial $u\in \Oqmat$ of degree $k>1$ can be written as
a product $u=u_1u_2$ with $\deg(u_{1,2})<k$
\begin{align*}
  x(u)=x(u_1)\vep(u_2)+\vep(u_1)x(u_2)=0
\end{align*}
as $x(u_1)=x(u_2)=0$ and therefore $x=0$.
\end{proof}
  
\begin{lemma}\label{prim}
$P(\Ubar )=\Lin \{F_{\beta _i},E_{\beta '_i}\,|\,\beta_i,\beta_i'\in
                                                            \Phi_r^+\}$.
\end{lemma}

\begin{proof}
By Lemma \ref{invK} and Lemma \ref{EFprim}
$P(\Ubar_+\oplus\Ubar_-)=P(\Ubar_+)\oplus P(\Ubar_-)$.
Suppose that $u=\sum \lambda _{n_1\dots n_M m_1\dots m_M}
\prod (F_{\beta_i})^{n_i} \prod (E_{\beta_i'})^{m_i} \in P(\Ubar )$,
$M=r(N-r)$ and $u\notin \Ubar_+\oplus\Ubar_-$.
By Lemma \ref{invK} one may assume that $u$ is homogeneous with respect to
the $\Z^{N-1}$-grading. Let $S_u\subset \N^{M}_0$ denote the subset
defined by
\begin{equation*}
  S_u:=\{(m_1,\dots,m_M)\,|\,\exists (n_1,\dots,n_M)\text{ such that }
                            \lambda_{n_1\dots n_M m_1\dots m_M}\neq 0\}.
\end{equation*}
Choose a multiindex $(k_1',\dots,k_M')\in S_u$ such that
  $\prod (E_{\beta_i'})^{k_i'}$ is maximal among the
  $\prod (E_{\beta_i'})^{l_i}$, $(l_1,\dots,l_M)\in S_u$ with respect to the
  filtration $\F_2$. By assumption
  $\prod (E_{\beta_i'})^{k_i'}\neq 1$. Pick $(k_1,\dots,k_M)$ such that
  $\lambda_{k_1\dots k_M k_1'\dots k_M'}\neq 0$.
  Write $\kopr u\in \Ubar\ot \Ubar$ with respect to the basis given in
  Proposition \ref{basprop} in the first tensor factor. The second
  tensor factor corresponding to $\prod (F_{\beta_i})^{k_i}$ is given by
  \begin{equation*}
    \sum_{(m_1,\dots,m_M)} \lambda _{k_1\dots k_M m_1\dots m_M}
     \prod (E_{\beta_i'})^{m_i}\neq 0
  \end{equation*}   
  as $u$ is homogeneous with respect to the $\Z^{N{-}1}$-grading.
  But this means that $u$ cannot be primitive.
\end{proof}

\begin{lemma}\label{nonzero}
  For any $x\in P(\Ubar)\setminus\{0\}$ the functional
  $\pair{\cdot}{x}:\B\rightarrow \C$ is nonzero.
\end{lemma}

\begin{proof}
By Lemma \ref{invK} and Lemma \ref{prim} the primitive elements of
$\Ubar$ form a direct sum of two non isomorphic irreducible right $K$
modules with highest weight vectors $E_r$ and $F_{\alpha_{1,N-1}}$.
Therefore one can find an element $k\in K$ such that $x k=E_r$ or
$x k= F_{\alpha_{1,N-1}}$. Restrict to the first case. The relation
\begin{align*}
\pair{b}{E_{r}}=\pair{b}{xk}=\pair{k\lact b}{x}
\end{align*}
implies that it suffices to find an element $b\in\B$ such that
$\pair{b}{E_{r}}\neq 0$. This is achieved by choosing $b=z_{r,r+1}$
as follows from (\ref{Eact}) and (\ref{zeps}).
\end{proof}  
  
\begin{lemma}\label{nondeg}
  The pairing $\pair{\cdot }{\cdot }: \B \ot \Ubar \to \comp $
  is nondegenerate.
\end{lemma}

\begin{proof}
Recall that the elements of $\K ^+\U $ vanish on $\B $. Hence the pairing
$\pair{\cdot }{\cdot }: \B \ot \Ubar \to \comp $
is just the restriction of the pairing between $\A $ and $\U $ to
$\B $ in the first component. Since the pairing between $\A $ and $\U $ is
nondegenerate, $\Ubar $ separates the elements of $\B $. On the other
hand, let $\mathcal{I}$ denote the subspace
\begin{align}
\mathcal{I}:=\{f\in \Ubar \,|\,\pair{b_{(1)}}{f}b_{(2)}=0\quad
\text{for all $b\in \B $}\}
\end{align}
of $\Ubar $. Clearly, $\mathcal{I}$ is a right $\U $ submodule of
$\Ubar $. Moreover, $\mathcal{I}$ is the kernel of the coalgebra map
$\Ubar \to \B ^\circ $ and hence a coideal.
Suppose that $\mathcal{I}\not=\{0\}$ and let $f\in \mathcal{I}$
be of minimal degree $k$ with respect to the coradical filtration.
Since $f\in \ker \vep $ and
\begin{align}
\kopr (f)-1\ot f-f\ot 1\in \Ubar _{k-1}\ot \Ubar _{k-1}\cap
\left(\mathcal{I}\ot \Ubar + \Ubar \ot \mathcal{I}\right),
\end{align}
by the minimality of $k$ we conclude that $f$ has to be primitive.
This is a contradiction to Lemma \ref{nonzero}.
\end{proof}

We conclude this section with an auxiliary lemma which gives an 
upper bound for the dimension of $(\B^+)^k/ (\B^+)^{k+1}$.  
\begin{lemma}\label{ubdim}
  $\dim(\B^+)^k/ (\B^+)^{k+1}\le \begin{pmatrix}
                                   2r(N-r)+k-1\\
                                   k
                                 \end{pmatrix}$.
\end{lemma}  
\begin{proof}
  First note that (\ref{proj}) and (\ref{zeps}) imply that
  $z_{ij}^+\in(\B^+)^2$ if $i,j\le r$ or $i,j>r$. Take for instance
  $i=j=N$ then
  \begin{align*}
     z_{NN}^+=&q^{2N+1}\sum_{i=1}^N q^{-2i}z_{Ni}z_{iN}-\vep(z_{NN})=
              \sum_{i<N}q^{2N+1-2i}z_{Ni}z_{iN}+q(z_{NN}^+)^2+2z_{NN}^+.
  \end{align*}
  Consider now $x\in(\B^+)^k/(\B^+)^{k+1}$. Then it follows that $x$ can be
  written as a linear combination of monomials of degree $k$ in the generators
  $z_{ij}\in M,$
  \begin{align*}
    M:=\{z_{mn}\,|\, m\le r,\,n>r \mbox{ or } m> r,\,n \le r  \}.
  \end{align*}
  The Equations (\ref{refl}) and (\ref{proj}) yield linear relations between
  these monomials. Define an $\N_0^{N-1}$-valued vector space filtration
  $\F$ of $(\B^+)^k/(\B^+)^{k+1}$ by
  \begin{align*}
    \deg\left(\prod_{m=1}^k z_{i_m j_m}\right)=
                    \sum_{m=1}^k\mathrm{e}_{|i_m-j_m|},\qquad z_{i_m j_m}\in M
  \end{align*}
  where $\N_0^{N-1}$ is ordered lexicographically
  $(\mathrm{e}_1>\mathrm{e}_2>\dots>\mathrm{e}_{N-1})$.
  To prove the lemma it suffices to show that for each pair
  $z_{kl},z_{ij}\in M$ there exists $c_{ijkl}\in \C\setminus \{0\}$
  such that
  \begin{align}
    z_{ij}z_{kl}-c_{ijkl}z_{kl} z_{ij}\in (\B^+)^3\label{comeps3}
  \end{align}
  up to terms of lower degree with respect to the filtration $\F$.
  All occuring cases are checked by direct computation. As an example we
  consider the case $i<l<j=k$.
  By \ref{relations}\ref{+-}.6
  \begin{align*}
    z_{ij}z_{jl}=& q^{-1}z_{jl}z_{ij}+q^{-1}\qd z_{jj}z_{il}
                   -q^{-1}\qd \sum_{m<j}q^{2j-2m}z_{im}z_{ml}\\
                =& q^{-1}z_{jl}z_{ij}+q^{-1}\qd z_{jj}^+z_{il}
                   +q^{-1}\qd\sum _{m>j}q^{2j-2m}z_{im}z_{ml}
                   +q^{-1}\qd z_{ij}z_{jl}
  \end{align*}
  which implies
  \begin{align*}
    z_{ij}z_{jl}= & q z_{jl}z_{ij}+q \qd z_{jj}^+ z_{il}
                   +q \qd\sum _{m>j}q^{2j-2m}z_{im}z_{ml}.
  \end{align*}
  As the last term of the right hand side is of lower degree with respect
  to $\F$ and $z_{jj}^+,z_{il}\in(\B^+)^2$ the obtained relation is of the
  desired form (\ref{comeps3}).
\end{proof}  

  Let $\Ubar_k$ denote the elements of degree $k$ in $\Ubar$ with respect to
  the filtration $\F_1$.
\begin{cor}\label{knondeg}  
  The pairing $\pair{\cdot }{\cdot }: \B/(\B^+)^{k+1} \ot \Ubar_k \to \comp $
  is nondegenerate.
\end{cor}
\begin{proof}
  By Lemma \ref{ubdim} and Proposition \ref{basprop}
  \begin{align*}
    \dim \B/(\B^+)^{k+1} \le \sum_{l=0}^k  \begin{pmatrix}
                                   2r(N-r)+l-1\\
                                   l
                                 \end{pmatrix}=\dim(\Ubar_k).
  \end{align*}                               
  On the other hand $\Ubar_k|_{(\B^+)^{k+1}}=0$. Using Lemma \ref{nondeg}
  one obtains that $\B/(\B^+)^{k+1}$ separates $\Ubar_k$ and
  hence $\dim(\Ubar_k)= \dim( \B/(\B^+)^{k+1})$. 
\end{proof}  

\begin{cor}
  The coradical filtration of $\Ubar$ coincides with $\F_1$.
\end{cor}  
\begin{proof}
  By Proposition 11.0.5 in \cite{b-Sweedler} and Lemma \ref{nondeg}
  one has
  \begin{equation*}
    C_k=\{f\in\Ubar\,|\, \pair{(\B^+)^{k+1}}{f}=0\}
  \end{equation*}  
  where $C_k$ denotes the elements of $\Ubar$ of degree $k$ with respect to
  the coradical filtration. By Corollary
  \ref{knondeg} the right hand side coincides with $\Ubar_k$.
\end{proof}

\section{Graded representations of $\oqgr$}\label{gradrep}

Let $V_0$ denote the subalgebra of $\B=\oqgr$ generated by the little
generators $z_{ii},i=1,\dots,N$, and let $V_+$ (resp. $V_-$)
denote the subalgebra of $\B^+$ generated by $z_{ij},i<j$
(resp. $i>j$).
\begin{lemma}\label{onto}
  The map $V_-\otimes V_0\otimes V_+\rightarrow \B$ given by multiplication
  is onto.
\end{lemma}
\begin{proof}
  Note first that the left hand side of the expressions in
  \ref{relations}.\ref{+0}, \ref{relations}.\ref{0-} and
  \ref{relations}.\ref{+-} consist of all possible products of
  generators $z_{db}z_{ca}$ with $z_{db}\in V_\alpha, z_{ca}\in V_{\beta}$
  and $(\alpha,\beta)\in\{(+,0),(0,-),(+,-)\}$.
  Let now $W_k$ denote the vector space generated by all $k$-fold products
  of the elements $z_{ij}, i,j=1,\dots,N$. It suffices to show that any
  element of $W_k$ can be written as a linear combination of monomials
  \begin{equation}
    v_{-,1}\dots v_{-,k_-} v_{0,1}\dots v_{0,k_0} v_{+,1}\dots v_{+,k_+}
                                                       \label{standard}
  \end{equation}
  in standard form where $v_{\alpha,j}$ denote generators of $V_\alpha$ for all
  $j=1,\dots,k_\alpha,\alpha\in\{-,0,+\}$ and $k=k_-+k_0+k_+$.
  This is achieved by induction over $k$ and the relation
  \begin{equation*}
    W_{k+1}=\sum_{c,a}W_k z_{ca}.
  \end{equation*}  
  Assume first $w\in W_k$ and $c=a$. By \ref{relations}.\ref{+0} one obtains
  \begin{align*}
    w z_{cc}\in W_k V_+
  \end{align*}  
  if $w\in \sum_{i<j}W_{k-1}z_{ij}$.

  \noindent On the other hand, if $c>a, c-a=l$ and $w\in W_k$ is a monomial
  (\ref{standard}) then \ref{relations}.\ref{0-} and \ref{relations}.\ref{+-}
  imply
  \begin{align*}
    w z_{ca}\in W_kV_0+W_kV_++\sum_{j<i<l+j}W_kz_{ij}.
  \end{align*}  
  Induction over $l$ yields the result.
\end{proof}  
 Let $T$ denote the quantum tangent space of a covariant first order
 differential calculus over $\B$, and let $\{t_i\,|\,i=1,\dots,d\}$ be a
 basis of
 $T^\vep=T\oplus \C\vep$. By Corollary \ref{corresp} one has 
 $\kow T^\vep\subset\B^\circ \ot T^\vep$.
 The elements $a_{ij}\in \B^\circ$ defined by
 \begin{equation*}
   \Delta t_i=a_{ij}\otimes t_j
 \end{equation*}  
 generate a finite dimensional subcoalgebra $T'\subset\B^\circ$.
 Thus they can be considered as matrix coefficients of a finite
 dimensional representation of $\B$.
\begin{lemma}\label{KinvT}
 Let $K'\subset \ug$ denote the commutative subalgebra generated by the
 elements $K_i,i=1,\dots,N-1$. Then $T'K'\subset T'$.
\end{lemma}
\begin{proof}
 By Corollary \ref{corresp} one has $T^\vep K\subset T^\vep$,
 further by (\ref{Kdef})
 $K_j\in K$ for all $j$. Therefore one can define 
 $\lambda_{ij}^l\in\C$ by $t_i K_j=\lambda_{ij}^lt_l$.
 As $K_j$ is invertible in $\ug$ the matrix
 $ \lambda_{ij}^l$ is invertible for fixed $j$.
 Then
 \begin{align*}
   \Delta(t_i K_j)=a_{ik}K_j\otimes t_k K_j
                  =a_{ik}K_j\otimes \lambda_{kj}^lt_l
                  =\lambda_{kj}^la_{ik}K_j\otimes t_l.
 \end{align*}                 
 On the other hand
 \begin{align*}
   \Delta(t_i K_j)=\Delta(\lambda_{ij}^mt_m)
                  =\lambda_{ij}^ma_{ml}\otimes t_l.
 \end{align*}                 
 Comparing coefficients of the linear independent elements
 $t_l$ one obtains
 \begin{align*}
   \lambda_{kj}^l a_{ik}K_j=\lambda_{ij}^m a_{ml}
 \end{align*}  
 for all $l,i,j$. The invertibiliy of
 $ \lambda_{ij}^l$ implies $a_{ik}K_j\in T'$.
\end{proof}
Let $P_+$ (resp. $P_-$) denote the set of monomials in the elements
$z_{ij},i<j$ (resp. $i>j$). 
\begin{cor}\label{pmnilp}
 The matrix coefficients $a_{ij}$ vanish on all but finitely many elements of
 $P_+$ and $P_-$.
\end{cor}
\begin{proof}
Note that for any $x\in \B$ and for any $K_k\in K'$
\begin{align*}
  (a_{ij}K_k)(x)=a_{ij}(K_k\triangleright x).
\end{align*}  
Assume that $a_{ij}(x)\neq 0$ for inifitely many elements $x\in P_+$.
Among the elements $x$ are common eigenvectors of the $K_k\triangleright$ with
infinitely many different eigenvalues. Therefore $\dim(a_{ij}K')=\infty$ in
contradiction to Lemma \ref{KinvT} and $\dim T'<\infty$. 
\end{proof}  
In particular the generators $z_{ij}$, $ i\neq j$, act as nilpotent operators
on the representation determined by the matrix coefficients $a_{ij}$.
Such representations are also annihilated by certain powers of $V_0^+$.
\begin{lemma}\label{nullnilp}
  Let $W$ denote a finite dimensional representation of $\B$ such that
  the generators $z_{ij}$, $ i\neq j$, act nilpotently.
  Then the commutative nonunital subring  $V_0^+\subset\B^+$ also acts
  nilpotently on $W$.
\end{lemma}
\begin{proof}
  The proof is performed in several steps.

  {\bf Step 1a:}
  Consider the following ordering on the generators $z_{ij},i<j$ of $V_+$:
  \begin{align*}
    z_{ij}<z_{kl}\Longleftrightarrow (i<k)\mbox{ or }(i=k\mbox{ and }j<l).
  \end{align*}  
  Let $v$ denote a common eigenvector of all $z_{ii}$, i.~e.
  \begin{align*}
    z_{nn}v=\lambda_n v\qquad \forall n=1\dots N
  \end{align*}  
  such that
  \begin{equation}
      z_{ij}v=0\qquad\forall z_{ij}<z_{kl}.\label{kleineristnull}
  \end{equation}    
  Then $w\colon= z_{kl} v$ satisfies
  \begin{enumerate}
    \item[1.] $z_{nn}w=\mu_n w$ where
      \begin{align}\label{mu+}
        \mu_n=\left\{ \begin{array} {lll}
                         q^{-2}\lambda_k&\mbox{if}&n=k,\\
                         \lambda_l+(1-q^{-2})\lambda_k&\mbox{if}&n=l,\\
                         \lambda_n&\mbox{else.}
                       \end{array} \right.
      \end{align}               
    \item[2.] $z_{ij}w=0\qquad\forall z_{ij}<z_{kl}.$
  \end{enumerate}
By assumption this implies the existence of a common eigenvector $v_+$ of all
$z_{ii},i=1,\dots, N$ such that $V_+v_+=0$. Such a common eigenvector will be
called a maximal eigenvector.

\textit{Proof of Step 1a:}
The value of $\mu_n$ can be computed by means of the list in the Appendix
\ref{relations}. Consider for example the case $n=l$. Then by the fourth
relation in \ref{relations}.\ref{+0} one has
\begin{eqnarray*}
  z_{ll}z_{kl} v&=&z_{kl} z_{ll}v+(1-q^{-2})\sum_{i<l}q^{2l-2i}z_{ki}z_{il}v\\
                &=& z_{kl} z_{ll}v+(1-q^{-2})q^{2l-2k}z_{kk}z_{kl}v
                   +(1-q^{-2})\sum_{k<i<l}q^{2l-2i}z_{ki}z_{il}v
\end{eqnarray*}
By assumption and by means of \ref{relations}.\ref{zdbzba} the last term is
simplified to
\begin{align*}
  -\hat{q}^2\sum_{k<i<l}q^{2l-2i}z_{kk}z_{kl}v.
\end{align*}  
Combining this with the result in the case $n=k$ one gets
\begin{eqnarray*}
  z_{ll}w&=&\lambda_l w+ (1-q^{-2})\lambda_k w.
\end{eqnarray*}  
The second property follows at once from the second and third relation of
\ref{relations}.\ref{++}

{\bf Step 1b:}
  In analogy to Step 1a consider the following ordering on the generators
  $z_{ij},i>j$ of $V_-$:
  \begin{align*}
    z_{ij}<z_{kl}\Longleftrightarrow (i>k)\mbox{ or }(i=k\mbox{ and }j>l).
  \end{align*}  
  Let $v$ denote a common eigenvector of all $z_{ii}$, i.~e.
  \begin{align*}
    z_{nn}v=\lambda_n v\qquad \forall n=1\dots N
  \end{align*}  
  such that
  \begin{equation}
      z_{ij}v=0\qquad\forall z_{ij}<z_{kl}.\label{kleineristnullMinus}
  \end{equation}    
  Then $w\colon= z_{kl} v$ satisfies
  \begin{enumerate}
    \item[1.] $z_{nn}w=\mu_n w$ where
      \begin{align}\label{mu-}
        \mu_n=\left\{ \begin{array} {lll}
                         q^2\lambda_k+(1-q^2)q^{2k-2N-1}&\mbox{if}&n=k,\\
                         \lambda_l+(1-q^2)\lambda_k+(q^2-1)q^{2k-2N-1}
                                                        &\mbox{if}&n=l,\\
                         \lambda_n&\mbox{else.}
                       \end{array} \right.
      \end{align}               
    \item[2.] $z_{ij}w=0\qquad\forall z_{ij}<z_{kl}.$
  \end{enumerate}
As for the case of $V_+$ this implies the existence of a common eigenvector
$v_-$ of all $z_{ii},i=1,\dots, N$ such that $V_-v_-=0$. Such a common
eigenvector will be called a minimal eigenvector.
  
\textit{ Proof of Step 1b:}
As above the value of $\mu_n$ can be computed by means of the list in
Appendix \ref{relations}. Consider again the case $n=l$. Then by the second
relation in \ref{relations}.\ref{0-} and
the projector property (\ref{proj}) one has
\begin{align*}
  z_{ll}z_{kl} v&=q^2z_{kl} z_{ll}v+(q^2-1)\left[q^{2l-2N-1}z_{kl}-\sum_{i\ge l}
                                   q^{2l-2i}z_{ki}z_{il}\right]v\\
                &=q^2 z_{kl} z_{ll}v+(q^2-1)q^{2l-2N-1}z_{kl}v
                   -(q^2-1)z_{kl}z_{ll}
                   -q\hat{q}\sum_{i> l}q^{2l-2i}z_{ki}z_{il}v.
\end{align*}
By assumption and by means of \ref{relations}.9.1 one obtains

\begin{align*}
  \sum_{i=l+1}^{k-1}q^{2l-2i}z_{ki}z_{il}v=
(1-q^{2(l-k+1)}) (q^{-2}z_{kk}w -q^{2k-2N-3}w).
\end{align*}
Combination with the result in the case $n=k$ yields the desired
expression.

 The second property follows at once from the second and third
relation of \ref{relations}.8.

\textbf{Step 2a:} Let $v_+$ denote a maximal eigenvector. Then
there exists a subset $\{i_1<i_2<\dots<i_s\}\subset \{1,\dots,N\}$ such that
for all $k$ the eigenvalue $\lambda_k$ of $z_{kk}$ is given by
\begin{align*}
\lambda_k=\begin{cases}q^{2l-2s-1}&\text{if $k=i_l$ for some $l$}\\
                       0&\text{else.}
          \end{cases}
\end{align*}

\textit{Proof of Step 2a:}
It follows from (\ref{proj}) and \ref{relations}.\ref{zabzba}.1 that
\begin{align*}
  z_{kk} v_+ = & q^{2N+1} \sum_{j\ge k}q^{-2j}z_{kj}z_{jk}v_+\\
             = & q^{2N+1-2k}z_{kk}^2v_+ + \qd\sum_{j>k}q^{2N+1-2j}
                 \left[q^{-1} z_{jj}z_{kk}-q z_{kk}^2
                       -\qd\sum_{i=k+1}^{j-1}z_{kk}z_{ii}\right]v_+.
\end{align*}
Thus $\lambda_k=0$ or $q\lambda_k+\qd \sum_{j>k}\lambda_j = 1$.
Therefore there exists a subset $\{i_1<i_2<\dots<i_n\}\subset \{1,\dots,N\}$
such that
\begin{align*}
\lambda_k=\begin{cases}q^{2l-2n-1}&\text{if $k=i_l$ for some $l$}\\
                       0&\text{else.}
          \end{cases}
\end{align*}
The relation (\ref{trace}) implies
$\sum_{k=1}^N\lambda_k=(1-q^{-2s})/(q-q^{-1})$ which
leads to $n=s$ as $q$ is not a root of unity.

\textbf{Step 2b:}
Let $v_-$ denote a minimal eigenvector. Then
there exists a subset $\{i_1<i_2<\dots<i_r\}\subset \{1,\dots,N\}$ such
that for all $k$ the eigenvalue $\lambda_k$ of $z_{kk}$ is given by
$\lambda_k=q^{2k-2N-1}-\lambda_k'$ where
\begin{align*}
\lambda_k'=\begin{cases}q^{2l-2N-1}&\text{if $k=i_l$ for some $l$}\\
                       0&\text{else.}
          \end{cases}
\end{align*}

\textit{Proof of Step 2b:}
It follows from (\ref{proj}) and \ref{relations}.\ref{zabzba}.2 that
\begin{align*}
  z_{kk} v_- = & \sum_{j\ge k}q^{2N+1-2j}z_{kj}z_{jk}v_-\\
             = & q^{2N+1-2k}\lambda_k^2 v_- + \qd\sum_{j<k}q^{2N+1-2j}
                 \Bigg[-q^{-1} \lambda_j\lambda_k +q^{2k-2j-1}\lambda_j^2\\
              &\phantom{ q^{2N+1-2k}\lambda_k^2 v_- + \qd\sum_{j<k}q^{2N+1-2j}
                 \Big[}    +\qd \sum_{i<j}\lambda_i\lambda_k
                  -q^{2k-2j}\qd\sum_{i<j}\lambda_i\lambda_j\Bigg]v_-.
\end{align*}
This implies
\begin{align*}
  q^2\lambda_k-\lambda_{k+1}=q^{2N-2k-1}(q^2\lambda_k-\lambda_{k+1})
       \left(\lambda_{k+1}+\lambda_k-q\qd\sum_{j<k}\lambda_j\right).
\end{align*}
Now the result follows in analogy to the proof of Step 2a.

\textbf{Step 3:} Using Steps 1 and 2 we will now prove the claim of the
Lemma. Let $v_-$ denote a minimal eigenvector with eigenvalues $\lambda_i$
and pick $k$ minimal such that $\lambda_k\neq 0$. Then by Step 2b we have
$k\le r+1$ and $\lambda_k=q^{2k-2N-1}$. Multiplication with suitable
elements of $V_+$ transforms $v_-$ into a maximal eigenvector $v_+$. By
(\ref{mu+})
\begin{align*}
  z_{kk}v_+ = q^{2k-2l-2N-1}v_+
\end{align*}  
where $l\ge 0$. Further by Step 2a 
\begin{align*}
  z_{kk}v_+ = q^{-2m+1}v_+
\end{align*}
where $m\le s$. Comparison of the exponents yields $k=r+1$ as $q$ is not a
root of unity. This implies in particular that
\begin{align}\label{minev}
z_{kk}v_-=\begin{cases}q^{2k-2N-1}v_-&\text{if $k\ge r+1$}\\
                       0&\text{else.}
          \end{cases}
\end{align}
Let now $v$ be an arbitrary common eigenvector with eigenvalues $\nu_k$
of $z_{kk}$, $k=1,\dots,N$. Multiplication with suitable elements of $V_-$
transforms $v$ into a minimal vector $v_-$ with eigenvalues $\lambda_k$.
By the above considerations $\lambda_k$ is given by (\ref{minev}).
Note that the tuple $(\lambda_1,\dots,\lambda_N)$ is invariant under the
invertible transformation (\ref{mu-}). This implies $\nu_k=\lambda_k$ and
therefore all common eigenvectors correspond to the same set of eigenvalues
independently of the representation.
Hence for any common eigenvector $v$ of $V_0$ one obtains
$(z_{ii}-\vep(z_{ii}))v=0$.
\end{proof}

\section{Classification}\label{classification}
\begin{lemma}\label{I=B+}
  Let $\mathcal{I}\subset\oqgr$ denote the ideal generated by
  \begin{align*}
    \{z_{kl}\,|\, k\le r \text{ or } l\le r\}.
  \end{align*}  
  Then $\mathcal{I}=\B^+$.
\end{lemma}
\begin{proof}
  By Proposition 2.3 in \cite{a-Kolb} $\mathcal{I}$ is equal to the kernel of
  the projection
  \begin{align*}
    \pi:\B=\oqgr\rightarrow\C
  \end{align*}
  induced by the surjective Hopf algebra homomorphism
  \begin{align*}
    \SUN&\rightarrow \cO_q(\textrm{SU}(N{-}r))\\
     u^i_j&\mapsto
         \begin{cases} \epsilon(u^i_j)&\textrm{if }i\le r
                       \textrm{ or } j\le r,\\
                       u^{i-r}_{j-r}&\textrm{else.}
          \end{cases}             
   \end{align*}  
    Thus $\mathcal{I}\subset\B^+$
  has codimension $1$ in $\B$ and therefore $\mathcal{I}=\B^+$.
\end{proof}  

\begin{lemma}\label{feps2k=0}
  Let $k\in\N$ and $f\in\B^\circ$ be a functional such that $f(x)=0$ for
  all monomials
  \begin{align}\label{standardx}
    x=x_{-,1}\dots x_{-,k_-} x_{0,1}^+\dots x_{0,k_0}^+ x_{+,1}\dots
    x_{+,k_+}\in
    V_-\otimes V_0^+\otimes V_+
  \end{align}
  where $k_-+k_0+k_+\ge k$ and $x_{\alpha,l}$ denote generators of $V_\alpha$.
  Then $f|_{(\B^+)^{2k}}=0$.
\end{lemma}

\begin{proof}
  Take $y\in(\B^+)^{2k}.$ By Lemma \ref{I=B+} the element $y$ can be
  written as a linear combination of monomials in the generators $z_{ij}$
  with at least $2k$ indices smaller than $r+1$. If we write these
  monomials in standard form (\ref{standardx})using the relations
  (\ref{refl}) the number of such indices will not decrease. Hence the
  appearing monomials have at least $k$ factors $z_{ij}=z_{ij}^+$ with
  $i\le r$ or $j\le r$. 
\end{proof}  

Recall that $\Ubar=\U/K^+\U$ can be considered as a subset of $\B^\circ$.

\begin{lemma}
  Let $T$ denote the quantum tangent space of a covariant first order
  differential calculus over $\B$. Then $T\subset \Ubar$.
\end{lemma}
\begin{proof}
  By Corollary \ref{corresp}  $T^\vep=T\oplus \C\vep\subset \B^\circ$ is a
  right $K$-module and a left $\B^\circ$-coideal. Consider the coalgebra $T'$
  defined above Lemma \ref{KinvT}. It is generated by functionals
  $a_{ij}\in\B^\circ$ which are matrix coefficients of a finite dimensional
  representation of $\B$. By Corollary \ref{pmnilp} and Lemma \ref{nullnilp}
  there exists $k\in\N$ such that $a_{ij}(x)=0$ for all monomials
  \begin{align*}
    x=x_{-,1}\dots x_{-,k_-} x_{0,1}^+\dots x_{0,k_0}^+ x_{+,1}\dots
    x_{+,k_+}\in
    V_-\otimes V_0^+\otimes V_+
  \end{align*}
  where $k_-+k_0+k_+\ge k$ and $x_{\alpha,l}$ denote generators of $V_\alpha$.
  Therefore by Lemma \ref{feps2k=0} one has $a_{ij}|_{(\B^+)^{2k}}=0$.
  Corollary $\ref{knondeg}$ implies $a_{ij}\in\Ubar_{2k-1}$. As
  $T^\vep\subset T'$ the claim of the lemma follows.
\end{proof}  

\begin{theorem}\label{mainTh}
  Any covariant first order differential calculus $\Gam$ over  $\oqgr$
  is uniquely determined by its quantum tangent space $T(\Gam)$. If
  $\dim\Gam\le 2r(N-r)$ then $T(\Gam)$ belongs to the
  following list.
  \begin{enumerate}
    \item For any $N,r$:
         $ \begin{aligned}[t]&T_0=\{0\},& &\dim\Gam(T_0)=0\\
                        &T_+=\Lin\{E_{\beta_i}\,|\,\beta_i\in\Phi^+_r\},&
                                     &\dim\Gam(T_+)=r(N-r)\\
                        &T_-=\Lin\{F_{\beta_i}\,|\,\beta_i\in\Phi^+_r\},&
                                     &\dim\Gam(T_-)=r(N-r)\\
                        &T=T_-\oplus T_+,& &\dim\Gam(T)=2r(N-r).            
          \end{aligned}$
    \item In addition
       \begin{align*}
          &\text{if $N{=}2,r{=}1$:} & &T_{1,+}=\Lin\{E_1,E_1^2\},&
                                  &T_{1,-}=\Lin\{F_1,F_1^2\},\\
          &\text{if $4{\le} N{\le} 7$, $2{\le} r{\le} N{-}2$:} &
                     & T_{2,+}=T_+\oplus V_+,&
                     & T_{2,-}=T_-\oplus V_-,
       \end{align*}
       where $V_{\pm}$ denotes the $K$-invariant
       $r(r{-}1)(N{-}r)(N{-}r{-}1)/4$-dimensional subspace of $(\Ubar_\pm)_2$.
   \end{enumerate}  
\end{theorem}

\begin{proof}
  The proof is performed in several steps.

  \textbf{Step 1:} If $T(\Gam)\cap\Ubar_\pm\neq\{0\}$ then $P(\Ubar_\pm)\subset
  T(\Gam)$.

  \textit{Proof of Step 1:} We prove the assertion in the case $\Ubar_+$.
  Choose an element $0{\neq} u\in T(\Gam)\cap \Ubar_+$ of degree $k$ with
  respect to the coradical filtration. Then
  \begin{align}\label{ucorad}
    \kopr u-1\otimes u- u\ot 1\in \sum_l (\Ubar_+)_l\ot(\Ubar_+)_{k{-}l}.
  \end{align}  
  Thus $u$ is primitive or there exists an element
  $0{\neq} u'\in T(\Gam)\cap \Ubar_+$ of degree ${<}k$. Therefore we can
  assume $u$ to be primitive.
  As in $\Ubar_+$ one has
  \begin{align*}
    E_{r}F_i=0,&&E_{r}K_j=
     q^{{-}2\delta_{j,r}{+}\delta_{j,r{+}1}{+}\delta_{j,r{-}1}}
     E_{r},&&i{\neq}r,
  \end{align*}
  the irreducible highest weight right $K$-module $E_rK{\subset}
  P(\Ubar_+)$ is $r(N{-}r)$-dimen\-sional.
  By Lemma \ref{EFprim} $\dim P(\Ubar_+){=}r(N{-}r)$ and therefore
  $E_rK=P(\Ubar_+)$. Now $u\in T(\Gam)$ and $T(\Gam)K\subset T(\Gam)$
  implies $P(\Ubar_+)\subset T(\Gam)$.

  Note that $T(\Gam)$ is $\Z^{N{-}1}$-graded as it is invariant under the
  action of  all $K_i$, $i=1,\dots,N-1$.
  
  \textbf{Step 2:}
  Assume that there exist $u_+\in\Ubar_+\cap T(\Gam)$ and
                          $u_-\in\Ubar_-\cap T(\Gam)$, $u_\pm\neq 0$ and
  $\dim T(\Gam)\le 2r(N{-}r)$. Then $T(\Gam)=T=P(\Ubar)$.

  \textit{Proof of Step 2:} By Step 1 $P(\Ubar_+)\subset T(\Gam)$ and
   $P(\Ubar_-)\subset T(\Gam)$.  Since $\dim P(\Ubar_\pm)=r(N{-}r)$ the
   assertion follows.

  \textbf{Step 3:}
  Suppose that $u=\sum \lambda _{n_1\dots n_M m_1\dots m_M}
  \prod (F_{\beta_i})^{n_i} \prod (E_{\beta_i'})^{m_i} \in T(\Gam)$ and
  $u\notin \Ubar_+\oplus\Ubar_-$. Then $\dim T(\Gam)> 2r(N-r)$.

  \textit{Proof of Step 3:}
  Without loss of generality we can assume $u$ to be homogeneous with
  respect to the $\Z^{N{-}1}$-grading.
  Let $S_u$ denote the subset of $\N_0^M$ defined by
  \begin{align*}
    S_u:=\{(m_1,\dots,m_M)\,|\,\exists (n_1,\dots,n_M)\text{ such that }
                            \lambda_{n_1\dots n_M m_1\dots m_M}\neq 0\}.
  \end{align*}                          
  Choose a multiindex $(k_1',\dots,k_M')\in S_u$ such that
  $\prod (E_{\beta_i'})^{k_i'}$ is maximal among the
  $\prod (E_{\beta_i'})^{l_i}$, $(l_1,\dots,l_M)\in S_u$ with respect to the
  filtration $\F_2$ from Section \ref{uk+u}. By assumption
  $\prod (E_{\beta_i'})^{k_i'}\neq 1$. Pick $(k_1,\dots,k_M)$ such that
  $\lambda_{k_1\dots k_M k_1'\dots k_M'}\neq 0$.
  Write $\kopr u\in \Ubar\ot T(\Gam)$ with respect to the basis given in
  Proposition \ref{basprop} in the first tensor factor. The second
  tensor factor corresponding to $\prod (F_{\beta_i})^{k_i}$ is given by
  \begin{align*}
    u_+:=\sum_{(m_1,\dots,m_M)} \lambda _{k_1\dots k_M m_1\dots m_M}
     \prod (E_{\beta_i'})^{m_i}\neq 0
  \end{align*}   
  as $u$ is homogeneous with respect to the $\Z^{N{-}1}$-grading.
  Therefore $u_+\in\Ubar_+\cap T(\Gam)\neq\{0\}$.
  Similarly one obtains that $\Ubar_-\cap T(\Gam)\neq\{0\}$.
  Now Step 2 and $u\notin \Ubar_+\oplus\Ubar_-$ imply the claim.
  
  \textbf{Step 4:} By Steps 2 and 3 it remains to consider the cases
  where $T(\Gam)\subset \Ubar_+$ or $T(\Gam)\subset \Ubar_-$.
  Consider the case $T(\Gam)\subset \Ubar_+$.
  Recall from the proof of Lemma \ref{EFprim} that the
  right $K$-module $(\Ubar_+)_k$ is dual to the left $K$-module
  $\Oqmat_k$ of homogeneous elements of degree $k$ in $\Oqmat$.

  By Step 1 one has $P(\Ubar_+)\subset T(\Gam)$.
  In what follows assume that there exists
  $u\in T(\Gam) \cap (\Ubar_+)_k$ for some $k\ge 2$.
  Then the coproduct of $u$ can be
  written as in (\ref{ucorad}).
  If no summands in $(\Ubar_+)_l\ot(\Ubar_+)_{k{-}l}$ occur for some
  $l\in\{1,\dots,k-1\}$ then
  $u(x)=0$ for all $x\in\Oqmat_k$. This is a contradiction to the duality
  between $\Ubar_+$ and $\Oqmat$. Thus for each $l=1,\dots,k$
  there exists a nonzero $u_l\in T(\Gam)\cap (\Ubar_+)_l$.
  
  If $r{=}1$ or $r{=}N{-}1$ then $\Oqmat_2$ is an irreducible $K$-module of
  dimension $N(N{-}1)/2$. Thus $\dim(\Ubar_+)_1+\dim(\Ubar_+)_2\le 2(N{-}1)$
  if and only if $N{=}2$. This proves the theorem if $r{=}1$ or
  $r{=}N{-}1$.

  If $2{\le}r{\le} N{-}2$ then $\Oqmat_2=V_1\oplus V_2$ is the direct
  sum of two irreducible $K$-modules of dimensions
  \begin{align*}
    \dim V_1=r(r{+}1)(N{-}r)(N{-}r{+}1)/4,\qquad
    \dim V_2=r(r{-}1)(N{-}r)(N{-}r{-}1)/4.
  \end{align*}  
  Since $\dim V_1>r(N{-}r)$ the component $V_1$ can not be a subspace of
  $T(\Gam)$. Moreover $\dim V_2\le r(N{-}r)$ if and only if
  $4{\le} N{\le} 6$ or
  $N{=}7,r{=}2$. In the case $N{=}6,r{=}3$ one has $\dim V_2= r(N{-}r)$.
  Let $V_+\subset (\Ubar_+)_2$ denote the $K$-submodule dual to $V_2$.
  Then if $u\in(\Ubar_+)_k$, $k{\ge} 2$ and $\dim T(\Gam) \le 2r(N{-}r)$
  one has $4{\le} N{\le} 7$ and $T_{2,+}=(\Ubar_+)_1\oplus V_+\subset T(\Gam)$.
  In particular if  $N{=}6,r{=}3$ then $T(\Gam)=T_{2,+}$.

  It remains to consider the cases $r{=}2$ and $r{=}N{-}2$.
  Since $\Oqmat_3=V_1'\oplus V_2'$ where $V_{1,2}'$ are irreducible
  $K$-modules of dimesions
  $$\dim V_1'=2(N{-}2)(N{-}1)N/3,\qquad \dim V_2'=2(N{-}2)(N{-}1)(N{-}3)/3$$
  one obtains 
  $$\dim V_{1,2}'+\dim(\Ubar_+)_1+\dim V_+ >4(N{-}2).$$
  Therefore $k\ge 3$ would imply $\dim T(\Gam)>2r(N{-}r)$ and hence
  $T(\Gam)=T_{2,+}$ .   
\end{proof}  

In \cite{a-Kolb} two differential calculi $\Gamma_+$ and $\Gamma_-$ of
dimension $0<\dim(\Gamma_\pm)\le r(N-r)$ were constructed (Prop.~3.1 in
 \cite{a-Kolb}).
Theorem \ref{mainTh} implies in particular that these caluli coincide with
$\Gamma(T_+)$ and $\Gamma(T_-)$.
Therefore the differential calculus corresponding to the quantum tangent
space $T=T_-\oplus T_+$ is isomorphic to $\gamgr=\Gamma_+\oplus\Gamma_-$.

\begin{appendix}
\section{Relations of $\oqgr$}\label{relations}
For notational convenience set $\qd=q-q^{-1}$.
\begin{enumerate}
\item $V_0$ is commutative
  \begin{itemize}
  \item $b=d<a=c$:\quad$z_{aa}z_{dd}=z_{dd}z_{aa}$
  \end{itemize}
\item $V_0\otimes V_+$\label{+0}
  \begin{itemize}
    \item $b=d<c<a$:\quad$z_{bb}z_{ca}=z_{ca}z_{bb}$
    \item $b=d=c<a$:\quad$z_{bb}z_{ba}=q^{-2}z_{ba}z_{bb}-q^{-1}\qd
                           \sum_{i<b}q^{2b-2i}z_{bi}z_{ia}$
    \item $c<b=d<a$:\quad \begin{minipage}[t]{.5\linewidth}
                          \begin{tabbing}
                            $z_{bb}z_{ca}$ \= $=z_{ca}z_{bb}
                           -q\qd z_{cb}z_{ba}-\qd^2
                           \sum_{i<b}q^{2b-2i}z_{ci}z_{ia}$\\
                          \> $=z_{ca}z_{bb}-\qd z_{ba}z_{cb}$
                         \end{tabbing}\end{minipage}
    \item $c<b=d=a$:\quad$z_{aa}z_{ca}=z_{ca}z_{aa}+q^{-1}\qd
                           \sum_{i<a}q^{2a-2i}z_{ci}z_{ia}$
    \item $c<a<b=d$:\quad$z_{bb}z_{ca}=z_{ca}z_{bb}$
  \end{itemize}
\item $V_0\otimes V_-$\label{0-}
  \begin{itemize}
    \item $b=d<a<c$:\quad$z_{bb}z_{ca}=z_{ca}z_{bb}$
    \item $b=d=a<c$:\quad$z_{bb}z_{cb}=q^{2}z_{cb}z_{bb}+q\qd\sum_{i<b}
                          q^{2b-2i}z_{ci}z_{ib}$
    \item $a<b=d<c$:\quad\begin{minipage}[t]{.5\linewidth}
                          \begin{tabbing}
               $z_{bb}z_{ca}$ \= $=z_{ca}z_{bb}+q\qd z_{cb}z_{ba}+
                          \qd^2\sum_{i<b}q^{2b-2i}z_{ci}z_{ia}$\\
                              \> $=z_{ca}z_{bb}+\qd z_{ba}z_{cb}$
                          \end{tabbing}\end{minipage}     
    \item $a<b=d=c$:\quad$z_{bb}z_{ba}=z_{ba}z_{bb}-q^{-1}\qd\sum_{i<b}
                          q^{2b-2i}z_{bi}z_{ia}$
    \item $a<c<b=d$:\quad$z_{bb}z_{ca}=z_{ca}z_{bb}$                    
  \end{itemize}
\item $V_+\otimes V_-$\label{+-}
  \begin{itemize}
    \item $d<b<a<c$ :\quad $z_{db}z_{ca}=z_{ca}z_{db}$
    \item $d<b=a<c$ :\quad $z_{da}z_{ca}=q z_{ca}z_{da}$
    \item $d<a<b<c$ :\quad $z_{db}z_{ca}=z_{ca}z_{db}+\qd
                            z_{cb}z_{da}$
    \item $d=a<b<c$ :\quad $z_{ab}z_{ca}=q z_{ca}z_{ab}+\qd
                            z_{cb}z_{aa}+
                            \qd\sum_{i<a}q^{2a-2i}z_{ci}z_{ib}$
    \item $a<d<b<c$ :\quad $z_{db}z_{ca}= z_{ca}z_{db}
                            +\qd z_{cb}z_{da}$
    \item $d<a<b=c$ :\quad\begin{minipage}[t]{.5\linewidth}
                          \begin{tabbing}
                $z_{db}z_{ba}=$ \= $q^{-1} z_{ba}z_{db}
                            +q^{-1}\qd z_{bb}z_{da}$\\
                               \>$-q^{-1}\qd\sum_{j<b}q^{2b-2j}z_{dj}z_{ja}$
                    \end{tabbing}\end{minipage}
    \item $d=a<b=c$ :\quad\begin{minipage}[t]{.5\linewidth}
                          \begin{tabbing}
                 $z_{ab}z_{ba}=$ \= $z_{ba}z_{ab}
                            +q^{-1}\qd z_{bb}z_{aa}
                            +q^{-1}\qd\sum_{j<a}q^{2a-2j}z_{bj}z_{jb}$\\
                          \>$-q^{-1}\qd\sum_{j<b}q^{2b-2j}z_{aj}z_{ja}$
                       \end{tabbing}\end{minipage}    
    \item $a<d<b=c$ :\quad\begin{minipage}[t]{.5\linewidth}
                          \begin{tabbing}
                    $z_{db}z_{ba}=$ \= $q^{-1} z_{ba}z_{db}
                            +q^{-1}\qd z_{da}z_{bb}$ \\
                         \> $-q^{-1}\qd\sum_{j<b}q^{2b-2j}z_{dj}z_{ja}$
                        \end{tabbing}\end{minipage}
    \item $d<a<c<b$ :\quad $z_{db}z_{ca}= z_{ca}z_{db}
                            +\qd z_{da}z_{cb}= z_{ca}z_{db}+\qd z_{cb}z_{da}$
    \item $d=a<c<b$ :\quad $z_{ab}z_{ca}=q z_{ca}z_{ab}
                            +\qd z_{aa}z_{cb}
                            +\qd\sum_{i<a}q^{2a-2i}z_{ci}z_{ib}$     
    \item $a<d<c<b$ :\quad $z_{db}z_{ca}= z_{ca}z_{db}
                            +\qd z_{da}z_{cb}$
    \item $a<d=c<b$ :\quad $z_{cb}z_{ca}=q z_{ca}z_{cb}$
    \item $a<c<d<b$ :\quad $z_{db}z_{ca}=z_{ca}z_{db}$                       
  \end{itemize}
\item $V_+\otimes V_+$\label{++}  
  \begin{itemize}
    \item $d<b=c<a$ :\quad $z_{db}z_{ba}=q^{-1}z_{ba}z_{db}
                            -q^{-1}\qd\sum_{i<b}q^{2b-2i}z_{di}z_{ia}$
    \item $d=c<b<a$ :\quad $z_{cb}z_{ca}=q^{-1}z_{ca}z_{cb}$
    \item $b>d<c<a$ :\quad\begin{minipage}[t]{.5\linewidth}
                          \begin{tabbing}
                      $q^{\delta_{bc}}z_{db}z_{ca}=$ \=
                              $q^{\delta_{ab}}z_{ca}z_{db}
                       -\delta_{bc}\qd\sum_{j<b}q^{2b-2j}z_{dj}z_{ja}$\\
                       \> $-(a<b)\qd z_{cb}z_{da}$
                       \end{tabbing}\end{minipage} 
  \end{itemize}
  where $(a<b)=1$ if $a<b$ and $(a<b)=0$ else. 
\item By induction from 5.1 \label{zdbzba}
  \begin{itemize}
    \item $d<b=c<a$ :\quad \begin{minipage}[t]{.5\linewidth}
                          \begin{tabbing}
                       $z_{db}z_{ba}=$ \= $q^{-1}z_{ba} z_{db}
                         -\hat{q}\sum_{i=d+1}^{b-1}z_{ia}z_{di}
                         -q\hat{q}z_{dd}z_{da}$\\
                       \>$-q\hat{q}\sum_{i=1}^{d-1}q^{2d-2i}z_{di}z_{ia}$
                      \end{tabbing}\end{minipage}   
  \end{itemize}
\item By induction from 4.7 \label{zabzba}
  \begin{itemize}
    \item $d=a<b=c$ :\quad \begin{minipage}[t]{.5\linewidth}
                          \begin{tabbing}
                         $z_{ab}$\=$z_{ba}= z_{ba}z_{ab}
                            +q^{-1}\hat{q}z_{aa}z_{bb}
                            -\hat{q}^2\sum_{i=a+1}^{b-1}z_{aa}z_{ii}$\\
                        \>$-q\hat{q}z_{aa}z_{aa}
                            -q\hat{q}\sum_{i=a+1}^{b-1}z_{ia}z_{ai}$\\
                       \> $+q^{-1}\hat{q}\sum_{j=1}^{a-1}q^{2a-2j}z_{bj}z_{jb}
                     +q\hat{q}\sum_{j=1}^{a-1}q^{2a-2j}z_{aj}z_{ja}$\\
                     \>$-\hat{q}^2\sum_{i=a+1}^{b-1}\sum_{j=1}^{a-1}q^{2a-2j}
                     z_{ij}z_{ji}$
                    \end{tabbing}\end{minipage}   
   \item $d=a<b=c$ :\quad \begin{minipage}[t]{.5\linewidth}
                          \begin{tabbing}
                        $z_{ba}$\=$z_{ab}= z_{ab}z_{ba}
                            -q^{-1}\hat{q}z_{aa}z_{bb}
                            +\hat{q}^2\sum_{j<a}z_{jj}z_{bb}$\\
                        \>$-q\hat{q}\sum_{j<a}z_{jb}z_{bj}
                       +q^{-1}\hat{q}\sum_{j=a}^{b-1}q^{2b-2j}z_{aj}z_{ja}$\\
                        \>$-q^{2b-2a}\hat{q}^2\sum_{j<a}z_{jj}z_{aa}
                     +q^{2b-2a+1}\hat{q}\sum_{j<a}z_{ja}z_{aj}$\\
                        \>$-\hat{q}^2\sum_{i=1}^{a-1}\sum_{j=a}^{b-1}q^{2b-2j}
                     z_{ij}z_{ji}$
                          \end{tabbing}\end{minipage}  
  \end{itemize}
\item $V_-\otimes V_-$ \label{Vmin} 
  \begin{itemize}
     \item $d>b=c>a$ :\quad \begin{minipage}[t]{.5\linewidth}
                          \begin{tabbing}
                    $z_{db}z_{ba}$ \= $=q^{-1}z_{ba}z_{db}-
                              \hat{q}q^{-1}\sum_{j<b}q^{2b-2j}z_{dj}z_{ja}$\\
                        \> $=qz_{ba}z_{db}+q\hat{q}\sum_{j>b}q^{2b-2j}
                              z_{dj}z_{ja}-q^{2b-2N}\hat{q}z_{da}$
                           \end{tabbing}\end{minipage}   
     \item $d=c>b>a$ :\quad $z_{db}z_{da}=q z_{da}z_{db}$
     \item $b<d>c>a$ :\quad\begin{minipage}[t]{.5\linewidth}
                          \begin{tabbing}
                        $q^{\delta_{bc}}z_{db}z_{ca}=$ \=
                                   $q^{-\delta_{ab}}z_{ca}z_{db}-
                          \delta_{bc}\qd\sum_{j<b}q^{2b-2j}z_{dj}z_{ja}$\\
                        \> $-(a>b)\qd z_{da}z_{cb}$
                        \end{tabbing}\end{minipage}    
  \end{itemize}
  where $(a>b)=1$ if $a>b$ and $(a>b)=0$ else.
\item By induction from \ref{Vmin}.1
  \begin{itemize}
     \item $d>b=c>a$ :\quad \begin{minipage}[t]{.5\linewidth}
                          \begin{tabbing}
                        $z_{db}z_{ba}=$ \= $q z_{ba}z_{db}-
                              \hat{q}q^{-1}q^{2d-2N-1}z_{da}
                              +\hat{q}\sum_{j>b}z_{ja}z_{dj}$\\
                              \> $+q^{-1}\hat{q}z_{dd}z_{da}
                              +q^{-1}\hat{q}\sum_{j>d}q^{2d-2j}z_{dj}z_{ja}$
                             \end{tabbing}\end{minipage}     
  \end{itemize}  
\end{enumerate}

\end{appendix}

\bibliographystyle{amsalpha}
\bibliography{litbank2}

\end{document}